\documentclass{article}
\usepackage[utf8]{inputenc}
\title{On Optimization of $\frac{1}{2}$-Approximation Path Cover Algorithm}
\author{Junyuan Lin, Guangpeng Ren}
\date{August 2020}

\usepackage{algorithm}  
\usepackage{algorithmicx} 
\usepackage{multicol}
\usepackage{csquotes}
\usepackage[utf8]{inputenc}
\usepackage[english]{babel}
\usepackage{algpseudocode}  
\usepackage[font=scriptsize,labelfont=bf]{caption}
\usepackage{blindtext}
\usepackage{color}
\usepackage{comment}
\usepackage{wrapfig}
\usepackage{subfig}
\usepackage{graphicx}
\usepackage{amsmath}
\newtheorem{theorem}{Theorem}
\newtheorem{lemma}[theorem]{Lemma}
\newtheorem{corollary}{Corollary}[theorem]

\algrenewcommand\algorithmicrequire{\textbf{Input:}}
\algrenewcommand\algorithmicensure{\textbf{Output:}}
\usepackage[
backend=biber,
style=numeric,
sorting=nyt
]{biblatex}
\newlength{\BiblioSpacing}
\renewenvironment{thebibliography}[1]{
\begin{oldthebibliography}{#1}
\setlength{\parskip}{\BiblioSpacing}
\setlength{\itemsep}{\BiblioSpacing}
}
{
\end{oldthebibliography}
}
\begin{document}
\maketitle
\begin{abstract}
In this paper, we propose a deterministic algorithm that approximates the optimal path cover on weighted undirected graphs. Based on the $\frac{1}{2}$-Approximation Path Cover Algorithm by Moran et al., we add a procedure to remove the redundant edges as the algorithm progresses. Our optimized algorithm not only significantly reduces the computation time but also maintains the theoretical guarantee of the original $\frac{1}{2}$-Approximation Path Cover Algorithm. To test the time complexity, we conduct numerical tests on graphs with various structures and random weights, from structured ring graphs to random graphs, such as Erdos-Renyi graphs. The tests demonstrate the effectiveness of our proposed algorithm on graphs, especially those with high degree nodes, and the advantages expand as the graph gets larger. Moreover, we also launch tests on various graphs/networks derived from a wide range of real-world problems to suggest the effectiveness and applicability of the proposed algorithm.
\end{abstract}
 

\section{Introduction}
In graph theory, the path cover refers to a set of paths that cover all vertices of a graph, where every vertex belongs to only one path. To find the optimal path cover, we would use the path cover algorithm to select a set of paths of minimum or maximum edge weights of the graph. Such algorithms have been widely implemented in many general graph theory problems, such as on cocomparability graphs, cographs, interval graphs, block graphs, and permutation graphs \cite{cocom, cograph, interval, block}. In applications, path cover also provides solutions to solve real-world problems. In \cite{ecommerce}, it provides a path analysis of site visiting to help understand the website traffic and improve marketing strategies. In \cite{dijkstra}, it shares a similar idea to the path cover algorithm by finding the path that has the minimum weight, such as the shortest distance or the shortest time on a map. In \cite{testing}, the path cover can provide the minimum number of test paths to cover different structural coverage criteria. 
\\\\
Finding an optimal path cover in any graph is an NP-complete problem. Therefore, we hope to obtain approximation algorithms to estimate an optimal path cover.  In \cite{pathcover}, Moran et al. introduced three fundamental approximation algorithms for the maximum weighted path cover of a graph.  $\frac{1}{2}$-Approximation Covering Algorithm is the first one in the paper. It is a greedy algorithm that works with the undirected graphs, and it can obtain at least $\frac{1}{2}$ of the optimal path cover weight on the graph. There are two other path cover approximation algorithms introduced in the \cite{pathcover} that are linear in programming. They both guarantee $\frac{2}{3}$ weight in undirected graphs and directed graphs. These two algorithms share a similar idea by finding a degree-constrained subgraph of the primary graph that has maximum weight. The degree of node in the subgraph is constrained to be less or equal to 2, and there are only paths and cycles left in the subgraph. Each cycle has at least three edges, and there is at least one edge in this cycle whose weight is less than $\frac{1}{3}$ of the total weight in this cycle. Thereby, these results give the $\frac{2}{3}$ theoretical bound of the algorithm.
\\\\
In recent years, there are also works done on approximation algorithms that take different approaches. In \cite{primal}, it uses primal-dual method and designs a $\frac{1}{2}$-approximation algorithm called $VCP_3$($Vertex$ $Cover$ $P_3$) where $P_3$ is a path with 3 vertices. The main idea of $VCP_3$ is to remove redundant vertices as the algorithm goes. Its time complexity takes $O(mn)$ time, where $m$ is the number of edges and $n$ is the number of vertices. In \cite{ksub}, it provides an algorithm specifically for $MWVCC_k$($Minimum$ $Weight$ $Connected$ $k-Subgraph$ $Cover$). This algorithm is an optimization for the algorithm in \cite{primal} where it only considers when $k=2$, yet the algorithm in \cite{ksub} proves $\frac{1}{k-1}$-approximation for $MWVCC_k$ problem. Its time complexity is $O(n^{2} \cdot m)$.
\\\\
In applications of path cover, most algorithm approaches were derived from approximation algorithms by Moran et al. Depending on the property of the problems, generalizations and optimizations are made to boost the efficiency. In \cite{auto}, it provides a generalized algorithm specifically working towards the flow graph, which is a presentation of all possible paths in a program. In \cite{lin}, the authors propose an optimization approach which is an adaptive multigrid method for linear systems with weighted graph Laplacians. Using the $\frac{1}{2}$-Approximation Covering Algorithm, the authors were able to approximate the level sets of the smooth error and further increase the accuracy and robustness of building multilevel structures. 
\\\\
As the data sets become larger and more complex nowadays, it calls for more robust path cover approximating algorithms to efficiently handle problems on these large graphs/networks. In this paper, we propose a deterministic algorithm based on the $\frac{1}{2}$-Approximation Covering Algorithm. Our algorithm inherits several advantages from the original one: 1) It is a greedy algorithm and straight-forward to implement. 2) Every step of the algorithm is deterministic, and there is no randomness involved. 3) There is a theoretical guarantee on time complexity. On top of these advantages, we take a further step to sequentially optimize it by removing the redundant edges as the algorithm progresses. Adding this step can significantly drop the computational time while still guaranteeing the $\frac{1}{2}$ weight lower bound. In section 2, we focus on $\frac{1}{2}$-Approximation Covering Algorithm by introducing some preliminaries and time complexity as well as the pseudocode and analysis of the algorithm. In section 3, we focus on the optimized algorithm by discussing our motivations, the pseudocode, the analysis of maintaining the theoretical guarantees as well as the process of removing edges. In section 4, we provide the numerical results and visualizations of both algorithms on Watts-Strogatz graph, a deterministic structure graph in the setting of this paper, and Erdos-Renyi graph, a random graph, as well as graphs from real-world problems to show the effectiveness and advantages of the optimized algorithm.
\section{$\frac{1}{2}$-Approximation Covering Algorithm}
We first introduce some preliminaries that help explain the $\frac{1}{2}$-Approximation Covering Algorithm, denoted as algorithm 1, and the complexity of the algorithm. Then we then show the pseudocode and review the proof of $\frac{1}{2}$ theoretical bound from \cite{pathcover}. 
\subsection{Preliminaries and complexity}
Consider an undirected weighted graph $G$, which contains vertex set $V$, edge set $E$, and associating weights $W$. An edge $e=\{u, v\}$ in $E$ contains two end points, and $w(e)$ is the weight of $e$ from $W$. We also define $N$ to be the number of vertices in $V$ and $M$ to be the number of edges in $E$. We denote $OPT$ as the optimal (maximal) path cover of $G$. $\texttt{cover}$ as the approximated path cover of $OPT$ that is generated by the approximating algorithms. Later on, as we analyze the lower bound of $W$ in $G$ that the algorithms can obtain, we denote $d_A$($u$) to be edges of $\texttt{cover}$ incident to $u$, $\epsilon(p)$ to be the set of edges in \texttt{cover}, $S(e)$ to be the set of edges incident to $e$, and $E_p=\{e=(u, v) | e \in \epsilon(p), d_A(u)=1\}$ \cite{pathcover}. In addition, we denote $H$ to be the number of edges in the path cover and $K$ to be the number of paths in the path cover.
\\\\
Algorithm 1 has a $\frac{1}{2}$ theoretical bound on the weight of $OPT$, which is where the name comes from. Its time complexity is $O(M \cdot LogM)$ which sources from the fact that initialization, the execution time of the loop, and the output are in $O(M)$, but the sorting time is in $O(M \cdot LogM)$. Thus, the entire test is eventually in $O(M \cdot LogM)$.
\subsection{Pseudocode}
Algorithm 1 is the fundamental path cover algorithm in \cite{pathcover}. Each step is deterministic and easy to follow. While the advantages are obvious, algorithm 1 can be time-consuming as the graph gets larger.
\begin{algorithm}[H]
  \caption{$\frac{1}{2}$ Approximation Path Cover}
  \begin{algorithmic}[1]
\Statex
    \Procedure{[\texttt{COVER}]}{\texttt{PathCover(\textbf{A})}}
    \Require $\texttt{A}$——an undirected positive weighted graph$\textit {G=(V,E,W)}$
    \Ensure \texttt{cover}——a path cover of graph $G$
    \State{\texttt{sorted\_edges} $\leftarrow$ \textbf{Sort} the edges in descending order based on \textit{W}}
    \For {$ \textit{e(u, v)} \in \texttt{sorted$\_$edges }$}
     \If{ neither $u$ nor $v$ is in any paths in  \texttt{cover}}
      \State{Add \{$u$, $v$\} as a new path in \texttt{cover}}
     \ElsIf{$u$ is the endpoint of a path in \texttt{cover} and $v$ is not in any paths}
      \State{Append \{$v$\} to \texttt{cover} \{path that contains $u$\}}
     \ElsIf{$v$ is the endpoint of a path in \texttt{cover} and $u$ is not in any paths}
      \State{Append \{$u$\} to \texttt{cover} \{path that contains $v$\}}
     \ElsIf{$u$ and $v$ are the endpoints of different paths in \texttt{cover}}
      \State{Merge two paths}
     \EndIf
    \EndFor
    \EndProcedure
  \end{algorithmic}
\end{algorithm}
\noindent We can observe that algorithm 1 is essentially a greedy algorithm by checking each edge based on the weight in descending order to build the path cover. Eventually, the weight of path cover is at least $\frac{1}{2}$ weight of $OPT$. More rigorous details will be discussed in Section 2.3.
\subsection{Algorithm 1 Analysis}
The goal of the analysis is to show the following theorem of algorithm 1:
\begin{theorem}
The weight of the path cover is at least $\frac{1}{2}$ of the weight of $OPT$
$$w(\texttt{cover}) \geq \frac{1}{2}\cdot w(OPT) $$
\end{theorem}
To review the $\frac{1}{2}$ lower bound proof in \cite{pathcover}, we first make the following classification:
\begin{flushleft}

(1) $E_1$ =\{$e$=($u$, $v$) $|$ $e$ $\in$ \texttt{cover} $\cap$  $OPT$. \} 

(2) $E_2$ =\{$e$=($u$, $v$) $|$ $e$ $\in$ $OPT$ $-$\texttt{cover}, $d_A$($u$)=$d_A$($v$)=1. \}

(3) $E_3$ =\{$e$=($u$, $v$) $|$ $e$ $\in$ $OPT$ $-$\texttt{cover}, max\{$d_A$($u$), $d_A$($v$)\}=2. \} 
\end{flushleft}
To generalize, we first consider $e$($u$, $v$) $\in$ $E_2$. It is not hard to state that $w$($e$) $\leq$ min($w$(\texttt{cover}) because, otherwise, the algorithm would have taken $e$ to \texttt{cover}. Since $S(e)$ is the set of edges incident to $e$. Then $S(e)$=($e_1$, $e_2)$) for each $e$ $\in$ $E_2$, and we can claim $w$($e$) $\leq$ $\frac{1}{2}$($w$($e_1$)+$w$($e_2$). 
\\\\
If we consider $e$($u$, $v$) $\in$ $E_3$ whose end points connect to at most two edges in \texttt{cover}, by similar argument, we let $S_1(e)$=($e_1$, $e_2)$) to be the set of edges incident to $u$ and $S_2(e)$=($e_3$, $e_4)$) to be the set of edges incident to $v$. Thus, we can still claim that $w$($e$) $\leq$ $\frac{1}{2}$($w$($e_1$)+$w$($e_2$) or $w$($e$) $\leq$ $\frac{1}{2}$($w$($e_3$)+$w$($e_4$)
\\\\
Based on the above statements, we conclude that

\begin{flushleft}

$w$($OPT$)=$w$($E_1$)+$w$($E_2$)+$w$($E_3$)$\leq$ $w_1$+$w_2$+$w_3$
\\

$$w_1 =  \sum_{e \in E_1}^{} w(e) $$ $$ w_2 = \sum_{\substack{e \in E_2 \\ S(e)=(e_1, e_2)}}^{}\frac{1}{2}(w(e_1)+w(e_2)) $$ $$ w_3 = \sum_{\substack{e \in E_3 \\ S(e)=(e_1, e_2)}}^{}\frac{1}{2}(w(e_1)+w(e_2))$$

\end{flushleft}
\noindent To show the contribution of each edge $\in \texttt{cover}$ to $W$, the lower bound of $\frac{1}{2}$ can be proved by discussing the following 4 cases.
\\\\
Assume that $e$ is an edge in $\epsilon(p)-E_p$ which means $d_A(u)\neq$1
\\\\
$\textbf{Case 1: }$ If $e\in OPT$, $w(e)$ contributes only once in $w_1$ since $e\in OPT \cap \texttt{cover}$ and at most twice in $w_3$ since $E_3$ requires $d_A(u)$ = 2, and it is multiplied by $\frac{1}{2}$ each time in $w_3$. Thus, $e$ contributes at most 2 $\cdot$ $w(e)$ to $W$.
\\\\
$\textbf{Case 2: }$ If $e\notin OPT$, $w_1$ and $w_2$ do not contribute since $e \notin OPT$ and $d_A(u)$ = 2, and $w_3$ contributes 2 $\cdot$ $w(e)$ as well since each end point $u$ and $v$ appears at most twice in $w_3$.
\\\\
Then assume that $e\in E_p$ which means $d_A(u)$=1
\\\\
$\textbf{Case 3: }$ If $e\in OPT$, $w(e)$ appears once in $w_1$. There might be at most one edge in $E_2$ incident to $e$, then $e$ appears at most once in $w_2$, and it is multiplied by $\frac{1}{2}$. $e$ also appears at most once in $w_3$ where $e$ may be one of the two edges incident to an edge in $E_3$, and it is multiplied by $\frac{1}{2}$. Thus, $e$ contributes at most 2 $\cdot$ $w(e)$ to $W$.
\\\\
$\textbf{Case 4: }$ If $e\notin OPT$, $w(e)$ does not appear in $w_1$. It appears once in $w_2$ since it may incident to an edge from $E_2$, and it is multiplied by $\frac{1}{2}$. It also appears at most once in $w_3$ since there could be at most two edges from $E_3$ incident to one side of $e$. Thus, $e$ contributes at most $\frac{3}{2}$ $\cdot$ $w(e)$ to $W$.
\\\\
Therefore, we can conclude that each edge from $\texttt{cover}$ contributes at most 2 $\cdot$ $w(e)$ to $W$. To take a further step, we can claim that it is guaranteed that at least $\frac{1}{2}$ of $W$ can be obtained by $\texttt{cover}$. Such relation can be presented as the following:
\\\\
$$w(OPT)=w(E_1)+w(E_2)+w(E_3)\leq W\leq 2\cdot w(\texttt{cover})$$
\\\\
To rewrite the expression, we can obtain:
\\\\
$$\frac{1}{2}\cdot w(OPT) \leq w(\texttt{cover})$$
\section{Optimized Algorithm}
As algorithm 1 is straight-forward and deterministic to build the path cover, we are motivated to develop an approach that can decrease the computational and preserve the theoretical results by removing redundant edges.
\subsection{Motivations}
\begin{figure}[htb]
\centering
\includegraphics[width=6cm]{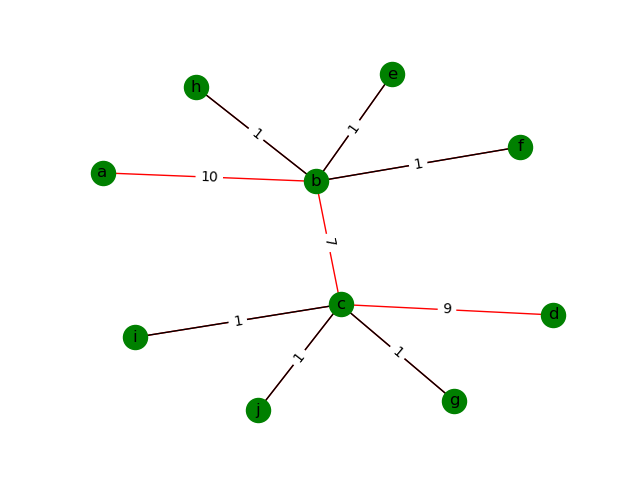}
\caption{Redundancy}
\label{fig:redundancy}
\end{figure}
\noindent In figure \ref{fig:redundancy}, edges in red are in path cover, and the edges in black are incident to the path. In algorithm 1, it would check all the edges to decide if they should be added to the path cover. However, by the definition of path, only edges incident to endpoints can be added to the path. Thus, checking these edges is the redundancy for the algorithm, and we add a step of removing these redundant edges in algorithm 2 to save computational time. The advantage of the time saved is significant as the average degree of the graph gets higher. Meanwhile, the time complexity and $\frac{1}{2}$ weight bound are still guaranteed in algorithm 2. 
\subsection{Pseudocode}
Steps 8, 11 and 14 of algorithm 2 are the steps that remove the redundant edges and save the time of checking these edges as the algorithm progresses. More numerical details on the efficiency increase will be discussed later in Section 4.
\begin{algorithm}[H]
  \caption{Sequential Optimized Path Cover}
  \begin{algorithmic}[1]
    \Statex
    \Procedure{[\texttt{COVER}]}{\texttt{PathCover(\textbf{A})}}
    \Require $\texttt{A}$——an undirected positive weighted graph$\textit {G=(V,E,W)}$
    \Ensure \texttt{cover}——a path cover of graph $G$
    \State{\texttt{sorted\_edges} $\leftarrow$ \textbf{Sort} the edges in descending order based on \textit{W}}
    \For {$ \textit{e(u, v)} \in \texttt{sorted$\_$edges }$}
     \If{ neither $u$ nor $v$ is in any paths in  \texttt{cover}}
      \State{Add \{$u$, $v$\} as a new path in \texttt{cover}}
     \ElsIf{$u$ is the endpoint of a path in \texttt{cover} and $v$ is not in any paths}
      \State{Append \{$v$\} to \texttt{cover} \{path that contains $u$\}}
      \State{Remove \{$adj(u), u$\} from \texttt{sorted\_edges}} \Comment{$adj(u)$ $\leftarrow$ adjacent nodes of $u$ in $G$ but not in \texttt{cover}}
     \ElsIf{$v$ is the endpoint of a path in \texttt{cover} and $u$ is not in any paths}
      \State{Append \{$u$\} to \texttt{cover} \{path that contains $v$\}}
      \State{Remove \{$adj(v), v$\}} from \texttt{sorted\_edges}\Comment{$adj(v)$ $\leftarrow$ adjacent nodes of $v$ in $G$ but not in \texttt{cover}}
     \ElsIf{$u$ and $v$ are the endpoints of different paths in \texttt{cover}}
      \State{Merge two paths}
      \State{Remove \{$adj(v), v$\} and \{$adj(u), u$\} from \texttt{sorted\_edges}}
     \EndIf
    \EndFor
    \EndProcedure
  \end{algorithmic}
\end{algorithm}
\noindent Comparing to the pseudocode of algorithm 1, these extra steps only remove the redundant edges from the loop of $\texttt{sorted$\_$edges }$ without influencing the build of path cover.
\subsection{Algorithm 2 Analysis}
As mentioned before, removing redundant edges does not affect the algorithm to obtain $\frac{1}{2}$ weight of $W$ since it does not interact with path cover in pseudocode. Therefore, the lower bound is still guaranteed in this new algorithm, as well as the time complexity.
\begin{figure}[H]
  \centering
  \begin{minipage}[]{0.4\textwidth}
    \includegraphics[width=\textwidth]{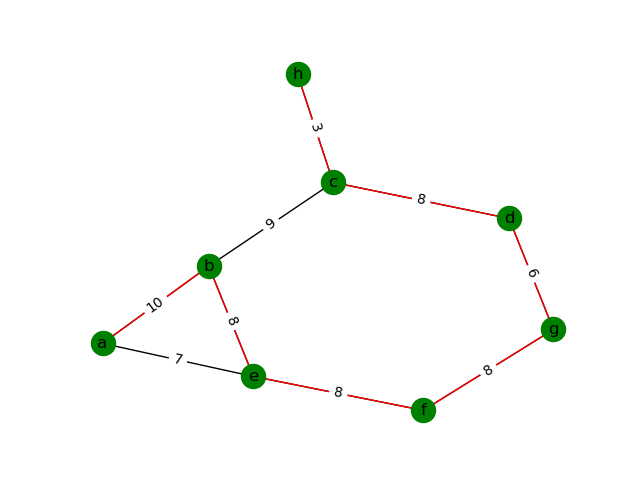}
    \caption{$OPT$}
    \label{fig:opt}
  \end{minipage}
  \hfill
  \begin{minipage}[]{0.4\textwidth}
    \includegraphics[width=\textwidth]{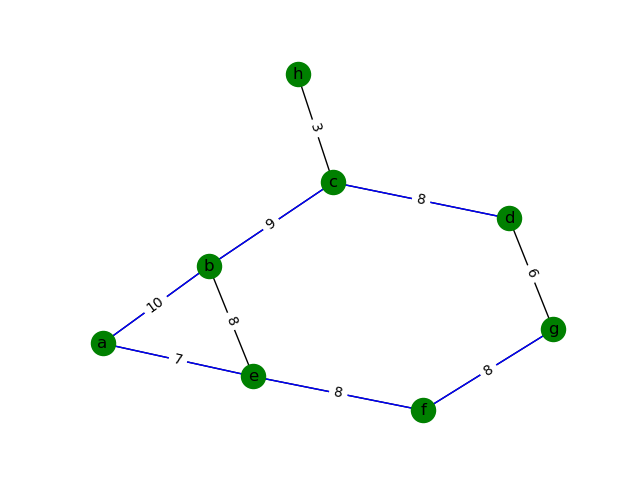}
    \caption{\texttt{cover} by Algorithm 1}
    \label{fig:eprime}
  \end{minipage}
\end{figure}
\noindent We generate an example of algorithm 2 following the same argument in the proof of algorithm 1. To visually compare, we generate graphs presented in figure \ref{fig:opt} and \ref{fig:eprime} which are the path covers generated by the optimal path cover($OPT$) which has the possible highest weight and the optimized algorithm, or \texttt{cover}. $OPT$ obtains the weight of 51, and \texttt{cover} obtains the weight of 50. The classification of edges can be made in the following:
\begin{flushleft}
(1) $E_1$ =\{e, f\}, \{f, g\}, \{c, d\}, \{a, b\}

(2) $E_2$ =\{d, g\}

(3) $E_3$ =\{b, e\} and \{c, h\}
\end{flushleft}
We first pick \{d, g\} from $E_2$ and compare it to the weight of each edge in $\texttt{cover}$, and we can claim that $w$(d, g) $\leq$ max($w(\texttt{cover}$)). Otherwise, the algorithm would have selected \{d, g\} instead of its two incident edges \{c, d\} and \{f, g\}. This can be expressed as $w$(d, g) $\leq$ $\frac{1}{2}$($w$(c, d)+$w$(f, g)). 
\\\\
Picking edges from $E_3$, we have \{b, e\} and \{c, h\}. By similar argument, we can claim the following expressions:
\\\\

$w$(b, e) $\leq$ $\frac{1}{2}$($w$(a, b)+$w$(b, c)) and $w$(c, h) $\leq$ $\frac{1}{2}$($w$(b, c)+$w$(c, d))
\\\\
Then we discuss the weight contribution of each edge. Looking back to the 4 cases mentioned previously, we can classify \{c, d\}, \{f, g\} for case 1, \{c, d\}, \{f, g\}, \{e, f\}, \{a, b\} for case 3 and \{a, e\}, \{b, c\} for case 4. 
\\\\
Edges in case 1 contribute once in $w_1$ since these edges are in $OPT$ and $\texttt{cover}$. Thus, these contributes at most 1 $\cdot$ $w(e)$ to $W$.
\\\\
Edges in case 3 contribute once in $w_1$ since they are all in $E_1$. There is only one edge \{d, g\} in $E_2$ incident to \{c, d\} and \{f, g\}. Thus, these two edges appear at most once in $w_2$, and their weight is multiplied by $\frac{1}{2}$. Edges that also appear at most once in $w_3$ as incident edges of $E_3$ are \{a, b\}, \{c, d\} and \{e, f\}, and they are multiplied by $\frac{1}{2}$ as well. Therefore, these edges in case 3 contribute at most 2 $\cdot$ $w(e)$ to $W$
\\\\
Edges in case 4 do not contribute in $w_1$ since they are not in $OPT$. They may appear once in $w_2$ if they incident to an edge from $E_2$, and the weight is multiplied by $\frac{1}{2}$. However, it does not apply in this example. Our selected edges in case 4 happened to be in $w_3$, so they contribute 1 $\cdot$ $w(e)$ which is within the $\frac{3}{2}$ bound mentioned earlier.
\\\\
This example still satisfies the expression we showed for algorithm 1, and it still works for algorithm 2:
\\\\
$$w(OPT)=w(E_1)+w(E_2)+w(E_3)\leq W\leq 2\cdot w(\texttt{cover})$$
\\\\
and
\\\\
$$\frac{1}{2}\cdot w(OPT) \leq w(\texttt{cover})$$
\\\\
Time complexity is still $O(M \cdot LogM)$. Previously we have initialization, the execution time of the loop, and the output is in $O(M)$, and the removal of redundant edges is also in $O(M)$. Also, the sorting time is still in $O(M \cdot LogM)$. Thus, the entire test is in $O(M \cdot LogM)$.
\subsection{Removed Edges}
To visualize the process of removing redundant edges in algorithm 2, we use the Watts-Strogatz graph in the setting of structured ring graphs. Thus, it is more straightforward to show the difference of both outputs from algorithms. 
\begin{figure}[H]
  \centering
  \begin{minipage}[]{0.3\textwidth}
    \includegraphics[width=\textwidth]{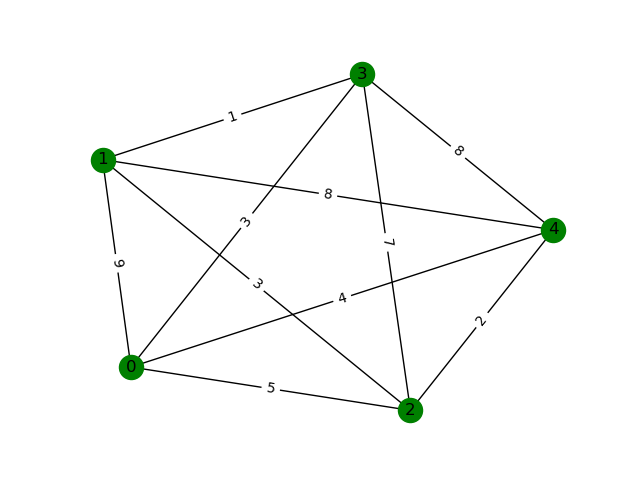}
    \caption{Initialized Watts-Strogatz Graph}
    \label{fig:iws}
  \end{minipage}
  \hfill
  \begin{minipage}[]{0.3\textwidth}
    \includegraphics[width=\textwidth]{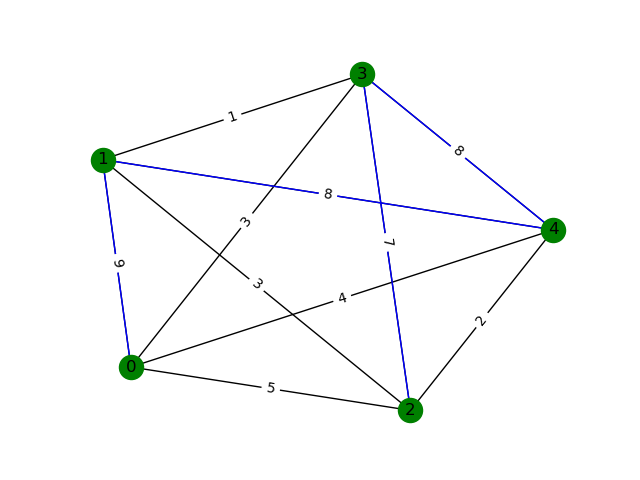}
    \caption{Output of Algorithm 1}
    \label{fig:ws1}
  \end{minipage}
  \hfill
  \begin{minipage}[]{0.3\textwidth}
    \includegraphics[width=\textwidth]{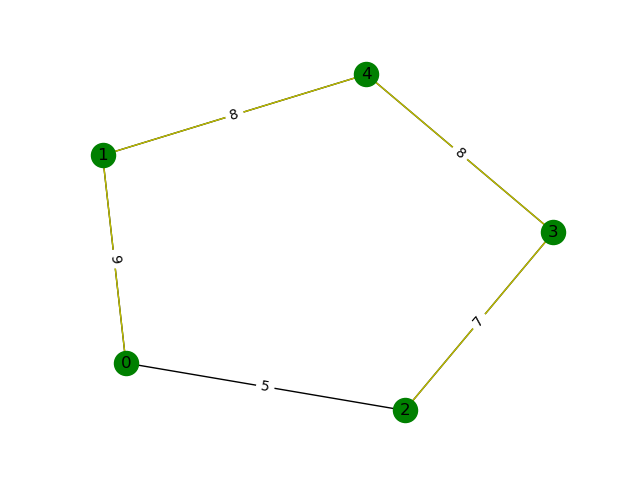}
    \caption{Output of Algorithm 2}
    \label{fig:ws2}
  \end{minipage}
\end{figure}
\noindent Figure \ref{fig:iws} is an initialized Watts-Strogatz graph with the degree of 4 to exemplify the process of generating path cover by two algorithms. Figure \ref{fig:ws1} is the output graph by algorithm 1, and the path cover is in blue. Figure \ref{fig:ws2} is the output graph by algorithm 2, and the path cover is in yellow.
\\\\
The process of algorithm 1 starts from figure \ref{fig:ws1} as the following. Based on the descending order of weights, it takes \{1, 0\} as the first edge, \{1, 4\} to connect node 1, \{3, 4\} to connect node 4, finally \{2, 3\} to connect node 3. Based on the weight, it should then check \{0, 2\}, but if this edge connects \{1, 0\} and \{2, 3\}, it is no longer a path. Therefore, algorithm 1 terminates.
\\\\
The process of algorithm 2 starting from figure \ref{fig:ws2} is mostly similar to algorithm 1. However, we need to remove redundant edges at the same time. It first takes \{1, 0\}, and there are no redundant edges yet since every other edge can be potentially chosen. Then it takes \{1, 4\}, then \{1, 3\} and \{1, 2\} cannot be taken because it needs to be a path. It takes \{3, 4\} next, then \{0, 4\} and \{2, 4\} cannot be taken for the same reason. Lastly, it removes \{0, 3\} after taking \{2, 3\}. \{0, 2\} is not the only remaining edge, and it cannot be taken. Therefore, algorithm 2 terminates.     
\\\\
We summarize our findings on the number of edges removed in the following claims:
\begin{lemma}\label{lm:1path}
In a connected weighted graph $G$, let $H$ be the number of edges in the path cover, $N$ be the number of vertices in the graph. If $G$ is covered by one single path, we have:
$$H = N-1$$
\end{lemma}
We can claim that if there is only one path in the path cover of the connected graph $G$, then $H$ should be $N-1$.
\\\\
According to this lemma, we can arrive at the following corollary:
\begin{corollary}
In a fully connected weighted graph $G$, denote $H$ as the number of edges in the path cover and $N$ as the number of vertices in the graph, we know that $$H = N-1$$.
\end{corollary}
This is because fully connected graphs are guaranteed to be covered by a single path. However, this is not guaranteed for any connected graphs.
\begin{figure}[H]
  \centering
  \begin{minipage}[]{0.45\textwidth}
    \includegraphics[width=\textwidth]{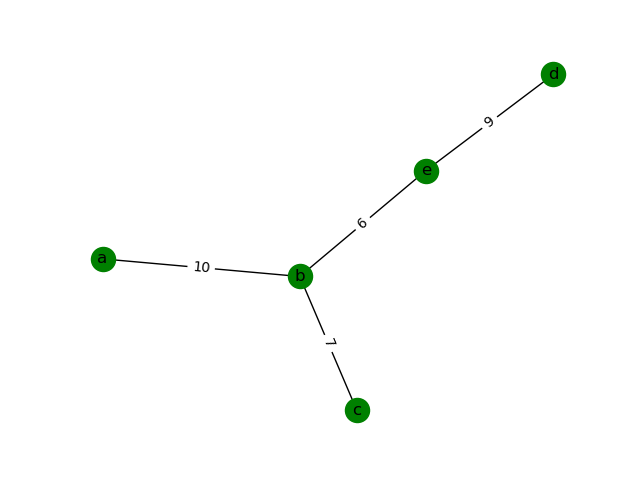}
    \caption{Original graph}
    \label{fig:toriginal}
  \end{minipage}
  \hfill
  \begin{minipage}[]{0.45\textwidth}
    \includegraphics[width=\textwidth]{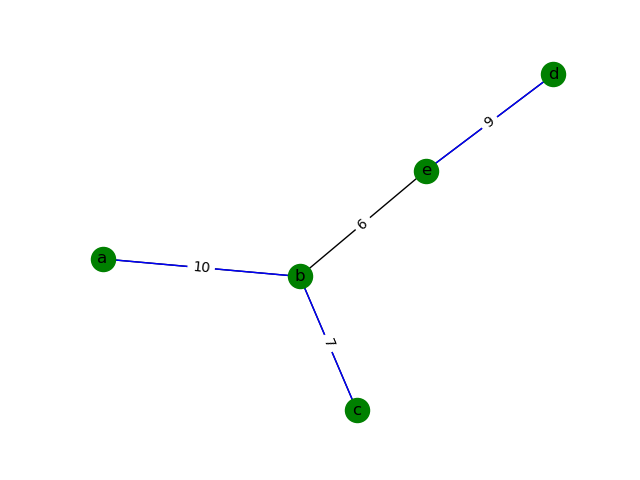}
    \caption{Algorithm 1(No edges removed)}
    \label{fig:t1}
  \end{minipage}
\end{figure}
\noindent For example, figure~\ref{fig:toriginal} is a connected graph but not a fully connected graph, we can see by launching algorithm 1 that there are two paths in the path cover (see figure~\ref{fig:t1}). In this case, $H$ does not equal $N-1$. For general cases where graphs are less connected, we use the following theorem to conclude the number of removed edges.
\begin{theorem}\label{thm:kpath}
In a weighted undirected graph, let $H$ be the number of edges in path cover, $N$ be the number of vertices in the graph, $K$ be the number of paths in path cover. $H$ can be expressed as the following:
$$H =N-K$$
\end{theorem}
\begin{figure}[H]
  \centering
  \begin{minipage}[]{0.45\textwidth}
    \includegraphics[width=\textwidth]{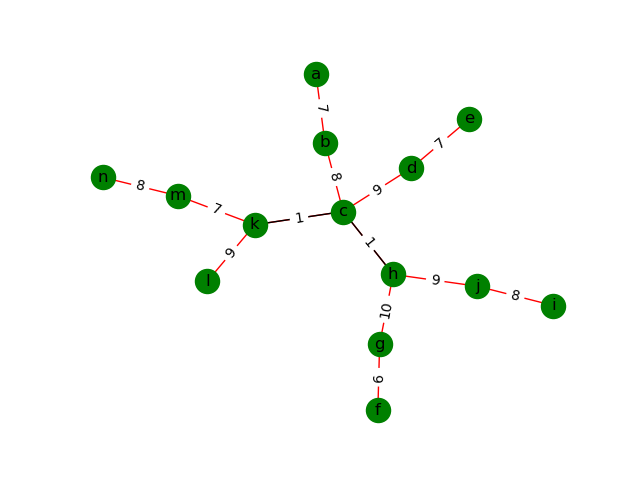}
    \caption{Theorem 3 visualization 1}
    \label{fig:thev1}
  \end{minipage}
  \hfill
  \begin{minipage}[]{0.45\textwidth}
    \includegraphics[width=\textwidth]{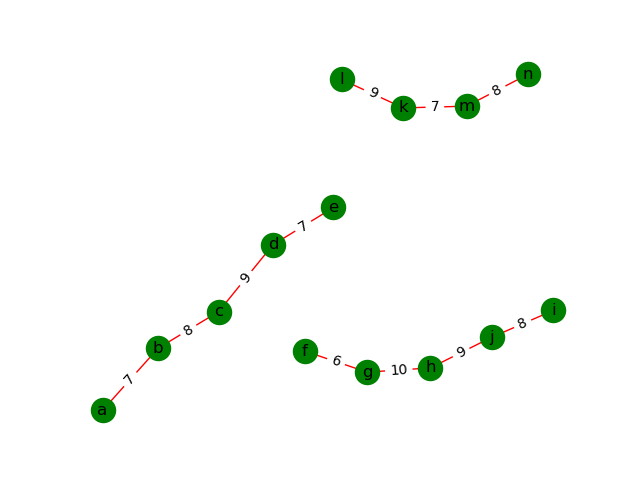}
    \caption{Theorem 3 visualization 2}
    \label{fig:thev2}
  \end{minipage}
\end{figure}
In figure \ref{fig:thev1}, it visualizes a graph containing three paths without removing the redundant edges, and figure \ref{fig:thev2} visualizes the graph that removes the redundant edges. To verify the theorem, we observe that $N=14$ and $K=3$. Thus $H=11$, and it verifies the formula given by Theorem 3.
\\\\
We include a direct proof for Theorem~\ref{thm:kpath}.
\\\\
$proof:$ To show $H=N-K$
\\\\
By Lemma~\ref{lm:1path}, 
$$H= \sum_{i=1}^{K} H_i=\sum_{i=1}^{K}(N_i-1)=\sum_{i=1}^{K}(N_i)-K=N-K$$
where $H_i$ is the number of edges and $N_i$ is the number of nodes in the $i-th$ path of path cover. Since all $i$ paths cover the entire vertex set, $\sum_{i=1}^{K}(N_i)=N$, Theorem 3 is therefore proved.

\section{Numerical Tests}
\begin{figure}[H]
  \centering
  \begin{minipage}[]{0.45\textwidth}
    \includegraphics[width=\textwidth]{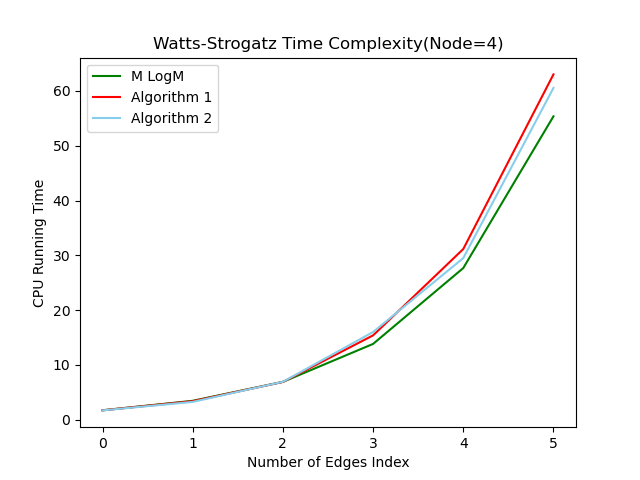}
    \caption{Watts-Strogatz graph Node=4}
    \label{fig:ws4}
  \end{minipage}
  \hfill
  \begin{minipage}[]{0.45\textwidth}
    \includegraphics[width=\textwidth]{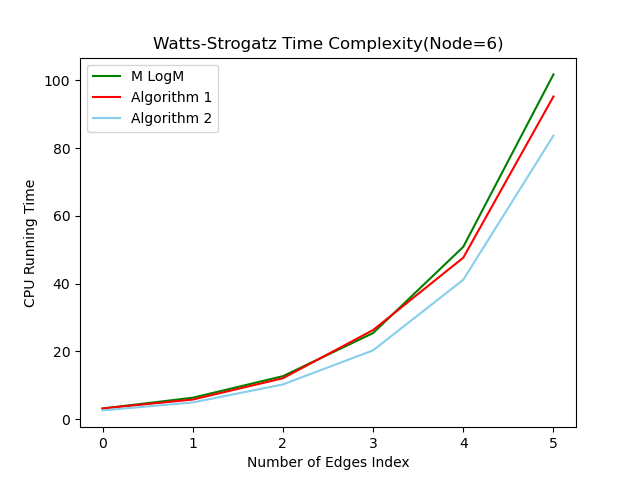}
    \caption{Watts-Strogatz graph Node=6}
    \label{fig:ws6}
  \end{minipage}
  \hfill
  \begin{minipage}[]{0.45\textwidth}
    \includegraphics[width=\textwidth]{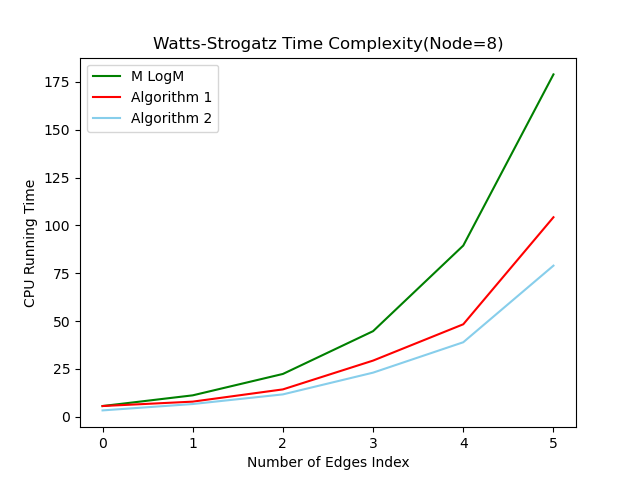}
    \caption{Watts-Strogatz graph Node=8}
    \label{fig:ws8}
  \end{minipage}
  \hfill
  \begin{minipage}[]{0.45\textwidth}
    \includegraphics[width=\textwidth]{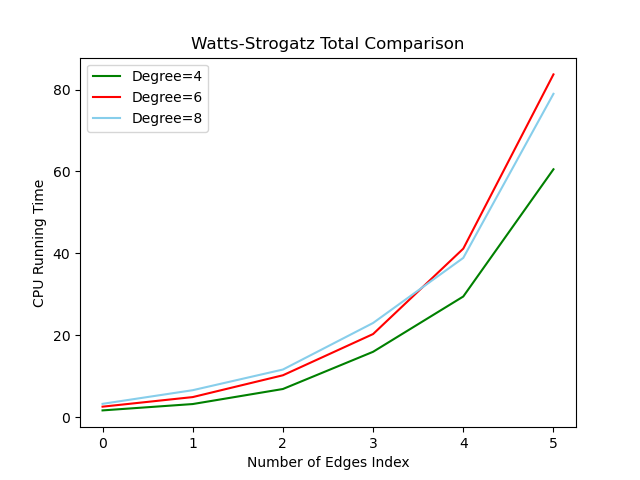}
    \caption{Watts-Strogatz Total Comparison}
    \label{fig:wst}
  \end{minipage}
\end{figure}

\begin{figure}[H]
  \centering
  \begin{minipage}[]{0.48\textwidth}
    \includegraphics[width=\textwidth]{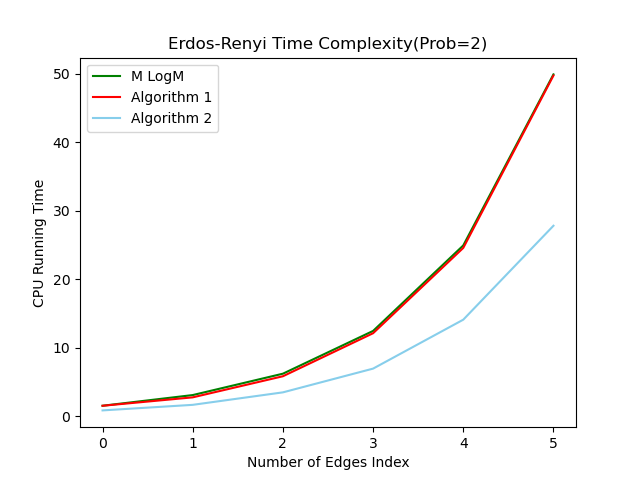}
    \caption{Erdos-Renyi graph Probability=2}
    \label{fig:er2}
  \end{minipage}
  \hfill
  \begin{minipage}[]{0.48\textwidth}
    \includegraphics[width=\textwidth]{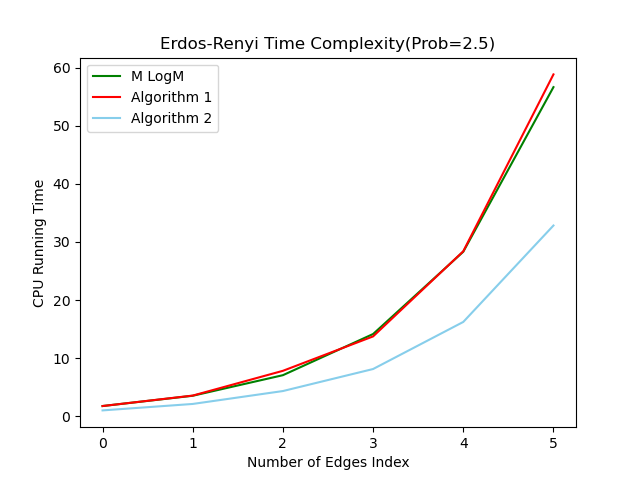}
    \caption{Erdos-Renyi graph Probability=2.5}
    \label{fig:er25}
  \end{minipage}
  \hfill
  \begin{minipage}[]{0.48\textwidth}
    \includegraphics[width=\textwidth]{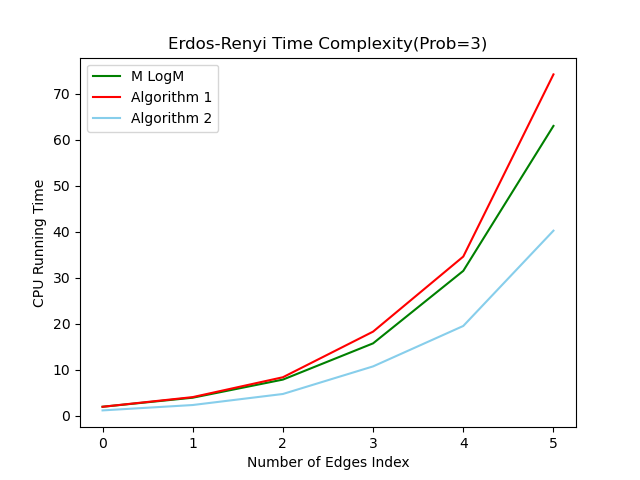}
    \caption{Erdos-Renyi graph Probability=3}
    \label{fig:er3}
  \end{minipage}
\end{figure}

\begin{figure}[H]
\centering
\subfloat{\includegraphics[width=5cm]{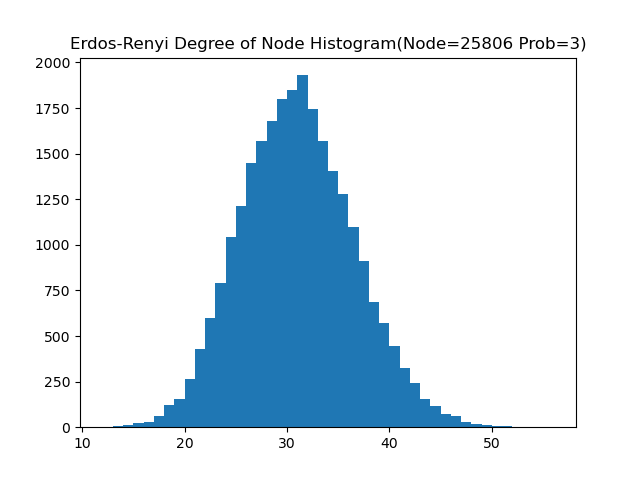}}\hfil
\subfloat{\includegraphics[width=5cm]{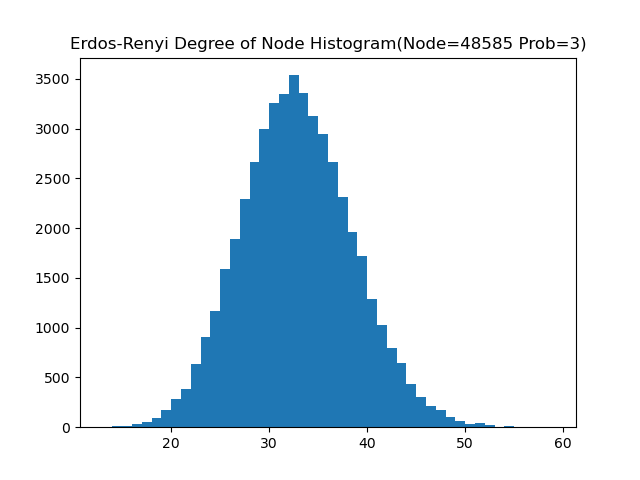}}\hfil 
\subfloat{\includegraphics[width=5cm]{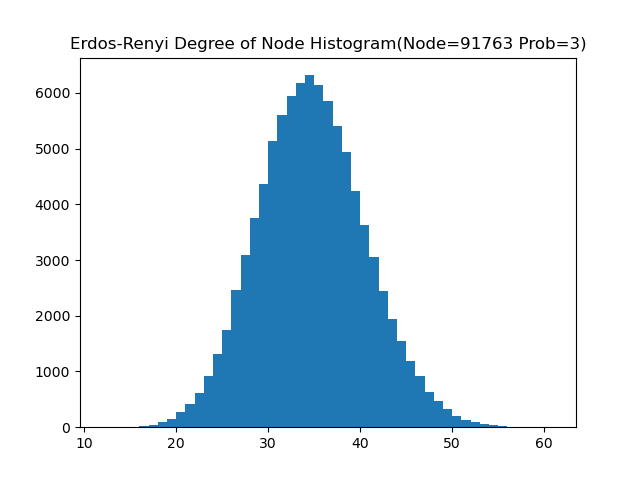}} \hfil
\subfloat{\includegraphics[width=5cm]{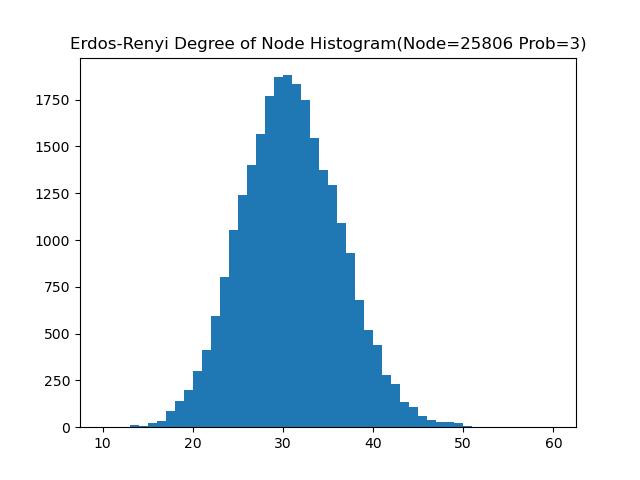}}\hfil   
\subfloat{\includegraphics[width=5cm]{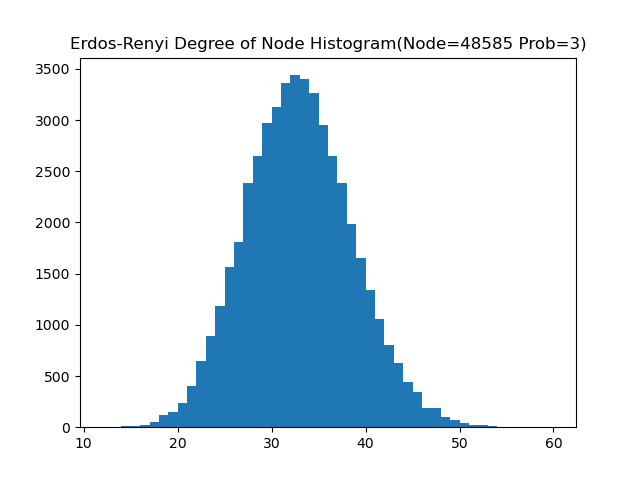}}\hfil
\subfloat{\includegraphics[width=5cm]{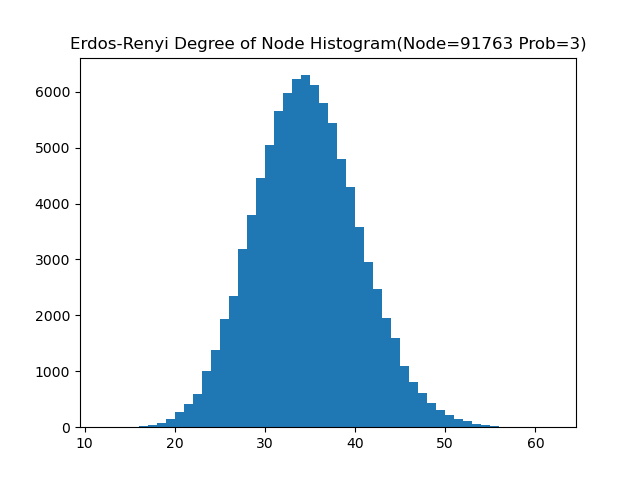}}\hfil
\caption{Least and Most Time Comparison}\label{ercom}
\end{figure}

The numerical tests are conducted with a 1.80 GHz Intel Core i7-8550U CPU, a quad-core processor, and 16 GB of RAM. The test program is written in Python where packages of networkx and numpy are imported. The dictionary is the main data structure throughout the program since it is implemented as hash tables, and its average time complexity is $O(1)$. At the end, packages of line-profiler and sys are imported to record the computational time.
\\\\
The following tables \ref{table:ws4}-\ref{table:er33} and figures \ref{fig:ws4}-\ref{ercom} provide the numerical results from the tests and visualizations on both Watts-Strogatz graphs and Erdos-Renyi graphs. We set 4, 6 and 8 as the degrees of node in Watts-Strogatz and $\frac{2 \cdot ln(M)}{M}$, $\frac{2.5 \cdot ln(M)}{M}$ and $\frac{3 \cdot ln(M)}{M}$, where $M$ is the number of edges, as the probability of generating edges in Erdos-Renyi is to see how it would affect the results as the density gets higher. 
\\\\

\begin{table}[H]
\caption{Watts-Strogatz Graph (Degree of nodes = 4)}
\centering
\begin{tabular}{llllll }  
\hline 
Nodes  & Edges   &  Avg. Degree & Algo 1 & Algo 2 & Time Ratio \\ \hline
131,072 & 262,144 & 4 & 1.73s & 1.71s & 0.99\\
262,144 & 524,288 & 4 & 3.39s & 3.26s & 0.96\\
524,288 & 1,048,576 & 4 & 6.90s & 6.93s & 1.00\\
1,048,576 & 2,097,152 & 4 & 15.40s & 16.01s & 1.04\\
2,097,152 & 4,194,304 & 4 & 31.15s & 29.50s & 0.95\\
4,194,304 & 8,388,608 & 4 & 63.01s & 60.54s & 0.96\\
\hline
\end{tabular}
\label{table:ws4} 
\end{table}
\begin{table}[H]
\caption{Watts-Strogatz Graph (Degree of nodes = 6)}
\centering
\begin{tabular}{llllll }  
\hline 
Nodes  & Edges   &  Avg. Degree & Algo 1 & Algo 2 & Time Ratio \\ \hline
131,072 & 393,216 & 6 & 3.18s & 2.62s & 0.82\\
262,144 & 786,432 & 6 & 5.82s & 4.96s & 0.85\\
524,288 & 1,572,864 & 6 & 12.11s & 10.28s & 0.85\\
1,048,576 & 3,145,728 & 6 & 26.33s & 20.34s & 0.77\\
2,097,152 & 6,291,456 & 6 & 47.68s & 41.17s & 0.86\\
4,194,304 & 12,582,912 & 6 & 95.24s & 83.69s & 0.88\\
\hline
\end{tabular}
\label{table:ws6} 
\end{table}
\begin{table}[H]
\caption{Watts-Strogatz Graph (Degree of nodes = 8)}
\centering
\begin{tabular}{llllll}  
\hline
Nodes  & Edges   & Avg. Degree & Algo 1 & Algo 2 & Time Ratio \\ \hline
131,072 & 524,288 & 8 & 5.59s & 3.31s & 0.59\\
262,144 & 1,048,576 & 8 & 7.89s & 6.64s & 0.84 \\
524,288 & 2,097,152 & 8& 14.29s & 11.68s & 0.82\\
1,048,576 & 4,194,304 & 8 & 29.34s & 23.04s & 0.79\\
2,097,152 & 8,388,608 & 8 & 48.31s & 38.92s & 0.81\\
4,194,304 & 16,777,216 & 8 & 104.24s & 78.95s & 0.76\\
\hline
\end{tabular}
\label{table:ws8} 
\end{table}
\noindent Tables \ref{table:ws4}-\ref{table:ws8} and figures \ref{fig:ws4}-\ref{fig:ws8} are the numerical results and visualizations of Watts-Strogatz graph. The setting of Watts-Strogatz in this paper is a deterministic structure graph as we let the degree of node to be fixed number 4, 6, and 8, and only the weights on edges are assigned randomly from 1 to 10. The starting point of $O(M \cdot LogM)$ is rescaled to the starting point of algorithm 1 to observe if algorithm 1 meets the time complexity bound, and algorithm 1 and algorithm 2 are plotted based on the values from the tables to compare the computational time. We can observe that algorithms 1 and 2 are slightly above the theoretical bound of $O(M \cdot LogM)$ when the degree of node is 4 due to the overhead of algorithms, but they have better performance as the degree of node gets higher. When the degree of node is 4, the graph is relatively sparse, and it is difficult to claim that algorithm 2 performs better than algorithm 1. However, as the degree of node increases, algorithm 2 starts to save more time, and the time ratio of time(algorithm 2) over time (algorithm 1) seems to be stable around 0.8. One notable result is that, in tables \ref{fig:ws6}-\ref{fig:ws8} and figure \ref{fig:wst}, time cost when the degree of nodes is 8 is less than when the degree of nodes is 6 which breaks the pattern that the computational time increases as degree of nodes increase starting from 2,097,152 nodes number. Therefore, in a deterministic graph, we can claim that algorithm 2 is expected to obtain better results as the graph and its average degree of nodes gets larger.
\begin{table}[H]
\caption{Erdos-Renyi Graph (Probability = $\frac{2\cdot(ln(M))}{M}$)}
\centering
\begin{tabular}{llllll}
\hline
Nodes  & Edges   &  Avg. Degree & Algo 1 & Algo 2 & Time Ratio \\ \hline
25,806 & 261,786 &  20.28  &  1.46s  & 0.92s & 0.63 \\
48,585 & 524,926 &  21.60  & 2.99s  &  1.84s  & 0.61  \\
91,763 & 1,048,862 & 22.86  &  6.22s  &  3.76s  &  0.60 \\
173,811 & 2,097,921 & 24.14 & 12.96s & 7.52s & 0.58 \\
330,076 & 4,195,049 & 25.42 & 23.19s & 14.10s & 0.61 \\ 
628,322 &  8,388,585 & 26.70 & 46.78s & 27.09s & 0.58 \\
 \hline
\end{tabular}
\label{table:er2} 
\end{table}
\begin{table}[H]
\caption{Erdos-Renyi Graph (Probability = $\frac{2.5\cdot(ln(M))}{M}$)}
\centering
\begin{tabular}{llllll}
\hline
Nodes  & Edges   &  Avg. Degree & Algo 1 & Algo 2 & Time Ratio \\ \hline
25,806 & 327,924 & 25.42 &  1.77s  & 1.02s &  0.58  \\
48,585 & 655,398 &  26.98  &  3.57s   & 2.13s  & 0.60  \\
91,763 & 1,310,775 & 28.56  &  7.83s  & 4.35s  & 0.56\\
173,811 & 2,623,118 & 30.62  & 13.75s & 8.14s  & 0.59\\
330,076 & 5,242,429 & 31.76  & 28.40s  & 16.24s  & 0.57\\ 
628,322 & 10,485,256  & 33.38  & 58.83s  & 32.82s  & 0.56  \\
\hline
\end{tabular}
\label{table:er25} 
\end{table}
\begin{table}[H]
\caption{Erdos-Renyi Graph (Probability = $\frac{3\cdot(ln(M))}{M}$)}
\centering
\begin{tabular}{llllll}
\hline
Nodes  & Edges   &  Avg. Degree & Algo 1 & Algo 2 & Time Ratio \\ \hline
25,806 & 393,057 & 30.46   & 1.97s   & 1.18s   & 0.60   \\
48,585 & 786,773 & 32.38   & 4.07s   & 2.35s   & 0.58  \\
91,763 & 1,572,307 & 34.26   & 8.39s   & 4.75s   & 0.57 \\
173,811 & 3,147,013 &  36.22  & 18.31s & 10.75s  & 0.59\\
330,076 & 6,288,240 & 38.10 & 34.63s & 19.54s & 0.56 \\ 
628,322 & 12,580,367 & 40.04 & 74.24s & 40.26s & 0.54\\
 \hline
\end{tabular}
\label{table:er3} 
\end{table}
\noindent While algorithm 2 is tested well on a deterministic structure graph, we would also like to see its performance on a random graph such as Erdos-Renyi graph. Erdos-Renyi graph generates its edges based on probability and assigns the weights on edges randomly from 1 to 10. In figures, \ref{fig:er2}-\ref{fig:er3}, algorithm 1 stays close to the theoretical bound, and algorithm 2 obviously has a slower increasing trend then algorithm 1 does. In tables \ref{table:er2}-\ref{table:er3}, the time ratio of algorithm 1 and 2 is stable in the range of 0.5 to 0.6 and likely to be decreasing. We can claim from the results that, as the graph and average degree get larger, algorithm 2 shows better performance.
\begin{table}[H]
s\caption{Erdos-Renyi Graph (Node=25806, Prob = $\frac{3\cdot(ln(M))}{M}$)}
\centering
\begin{tabular}{lllllll}
\hline
Edges  &  Avg. Degree & Skew. & Kurt. & Algo 1 & Algo 2 & Time Ratio\\ \hline
394,137& 30.54	&0.1940	&0.0185 &1.91s&  1.10s& 0.58\\
393,468& 30.50	&0.1961	&0.0075 &1.92s&  1.11s& 0.58\\
393,088 &30.46	&0.1876	&0.0639 &1.93s&  1.13s& 0.59\\
393,690& 30.52	&0.1647	&0.0707 &1.92s&  1.15s& 0.60\\
392,296& 30.40	&0.1783	&0.0148 &1.96s&  1.15s& 0.59\\
391,919&30.36  &0.1934 &0.0669 &2.05s&	1.21s& 0.59\\
393,543& 30.50	&0.1729	&0.0541 &1.99s&  1.21s& 0.61\\
393,203 &30.48	&0.1892	&0.0591 &1.98s&  1.22s& 0.62\\
392,256 &30.40	&0.1675	&0.0786 &2.07s&  1.27s& 0.61\\
392,968 &30.46	&0.1661	&0.0414 &1.97s&  1.27s& 0.64\\
 \hline
\end{tabular}
\label{table:er31} 
\end{table}
\begin{table}[H]
\caption{Erdos-Renyi Graph (Node=48585, Prob = $\frac{3\cdot(ln(M))}{M}$)}
\centering
\begin{tabular}{lllllll}
\hline
Edges  &  Avg. Degree & Skew. & Kurt. & Algo 1 & Algo 2 & Time Ratio\\ \hline
786571& 32.38	&0.1771	&0.0490 &4.10s&  2.19s& 0.53\\
787096& 32.40   &0.1827 &0.0238 &3.98s&	2.21s& 0.59\\
786128& 32.36	&0.1808	&0.0142 &4.11s&  2.24s& 0.55\\
786808& 32.38	&0.1649	&0.0602 &4.05s&  2.26s& 0.56\\
788167& 32.44	&0.1722	&0.0588 &4.08s&  2.31s& 0.57\\
784711& 23.30   &0.1982	&0.1131 &4.02s&  2.34s& 0.58\\
785326& 32.32	&0.1817	&0.0457 &4.20s&  2.38s& 0.57\\
788140& 32.44	&0.1774	&0.0486 &4.06s&  2.40s& 0.59\\
786928& 32.40	&0.1762	&0.0236 &4.07s&  2.49s& 0.61\\
787851& 32.44	&0.1746	&0.3380 &4.07s&  2.53s& 0.62\\
 \hline
\end{tabular}
\label{table:er32} 
\end{table}
\begin{table}[H]
\caption{Erdos-Renyi Graph (Node=91763, Prob = $\frac{3\cdot(ln(M))}{M}$)}
\centering
\begin{tabular}{lllllll}
\hline
Edges  & Avg. Degree & Skew. & Kurt. & Algo 1 & Algo 2 & Time Ratio\\ \hline
1571521& 34.26	&0.1649	&0.0092 &7.73s&  4.34s& 0.56\\
1571816& 34.26	&0.1737	&0.0157 &8.20s&  4.56s& 0.56\\
1573021& 34.28  &0.1618 &0.0542 &8.10s&	 4.61s& 0.57\\
1573531& 34.30	&0.1606	&0.0046 &8.23s&  4.63s& 0.56\\
1573614& 34.30  &0.1585	&0.0194 &8.56s&  4.67s& 0.55\\
1572042& 34.26	&0.1684	&0.0304 &8.15s&  4.77s& 0.59\\
1573326& 34.30	&0.1755	&0.0283 &8.16s&  4.87s& 0.60\\
1572499& 34.28	&0.1711	&0.0136 &8.29s&  4.91s& 0.59\\
1573169& 34.28	&0.1620	&-0.0068&9.25s&  5.01s& 0.54\\
1568528& 34.18	&0.1752	&0.0399 &9.19s&  5.16s& 0.56\\
 \hline
\end{tabular}
\label{table:er33} 
\end{table}
\noindent In a random graph, the variance of different test results while the nodes are the same could be high. To investigate what might affect the computational time of each test, we observe that the skewness and kurtosis of the degree of node distribution could be the reasons. In tables \ref{table:er31}-\ref{table:er33}, we explicitly list out 10 test results of three different nodes, and the list is sorted by the algorithm 2 time in ascending order. The idea of algorithm 2 is to eliminate the extra edges to save time, so we are looking for the graph whose nodes have a more high average degree of nodes. Hypothetically, we would prefer the degree of nodes distribution that has low skewness and kurtosis. In other words, we want the distribution skewing to the right with the heavy tail which means more nodes with a high degree. Looking at the tables, we find this to be mostly true while the density and edges number are close. In figure \ref{ercom}, the top 3 graphs are the distribution of degree from the least running time graph of three nodes number, and the bottom 3 graphs are the distribution of degree from the most running time graph of three nodes number. We can tell that the top 3 distributions skew more to the right and have heavier tails comparing to the bottom 3 accordingly.
\begin{table}[H]
\caption{gemsec-Facebook}
\centering
\begin{tabular}{llllllllll}
\hline
Group & Nodes  & Edges  &  Avg. Degree& Skew. & Kurt. & Algo 1 & Algo 2 & Time Ratio \\ \hline
Government & 7,057 & 89,429  & 25.34 & 6.02 & 64.63 & 0.45s & 0.26s & 0.58 \\
New Sites & 27,917 & 205,964  & 14.76& 8.83 & 154.37 & 1.18s & 0.66s  & 0.56 \\
Athletes & 13,866 & 86,811  & 12.52& 6.54 & 85.59 & 0.44s & 0.28s &  0.64\\
Public Figures & 11,565 & 67,038 &11.59 & 5.40  & 42.56 & 0.36s & 0.21s  &  0.58\\
TV Shows & 3,892 & 17,239 & 8.86& 3.71 & 18.21  & 0.10s & 0.06s &  0.60\\ 
Politician & 5,908 & 41,706 & 14.12& 4.46 & 33.79 & 0.28s  & 0.15s &  0.54\\
Artist & 50,515 & 819,090  &32.43 & 7.40 &  90.50 & 4.37s & 2.43s  &  0.56\\
Company & 14,113 & 52,126 & 7.39& 6.43 & 77.74 & 0.28s & 0.19s  &  0.68\\
 \hline
\end{tabular}
\label{table:fb} 
\end{table}
\begin{table}[H]
\caption{gemsec-Deezer}
\centering
\begin{tabular}{llllllllll}
\hline
Group & Nodes  & Edges  &  Avg. Degree& Skew. & Kurt. & Algo 1 & Algo 2 & Time Ratio \\ \hline
HR & 54,573 & 498,202  & 18.26& 2.91&21.63 &2.87s &1.66s & 0.58\\
HU & 47,538 & 222,887  & 9.38& 2.07& 9.03&1.26s &0.86s & 0.68\\
RO & 41,773 & 125,826  & 6.02& 3.27&23.64 & 0.72s &0.53s & 0.74\\
 \hline
\end{tabular}
\label{table:dee} 
\end{table}
\begin{table}[H]
\caption{Road Network}
\begin{tabular}{llllllllll}
\hline
Group & Nodes  & Edges  &Avg. Degree& Skew. & Kurt. & Algo 1 & Algo 2 & Time Ratio \\ \hline
TX & 1,379,917 & 1,921,660 & 2.79 & -0.59 &-0.48 & 15.04s& 14.97s& 1.00\\
CA & 1,965,206 & 2,766,607 & 2.82 & -0.55& -0.39&21.51s &21.49s & 1.00\\
PA & 1,088,092 & 1,541,898 & 2.83 &-0.58&-0.42 & 11.85s &12.19s &1.03 \\
 \hline
\end{tabular}
\label{table:road} 
\end{table}
\begin{figure}
\centering
\subfloat{\includegraphics[width=4cm]{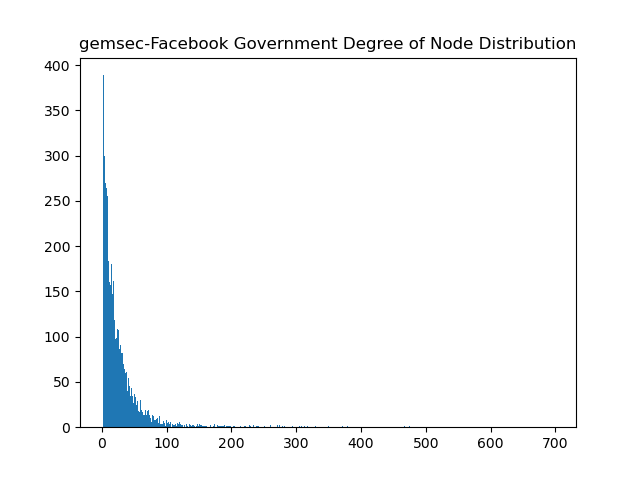}}\hfil
\subfloat{\includegraphics[width=4cm]{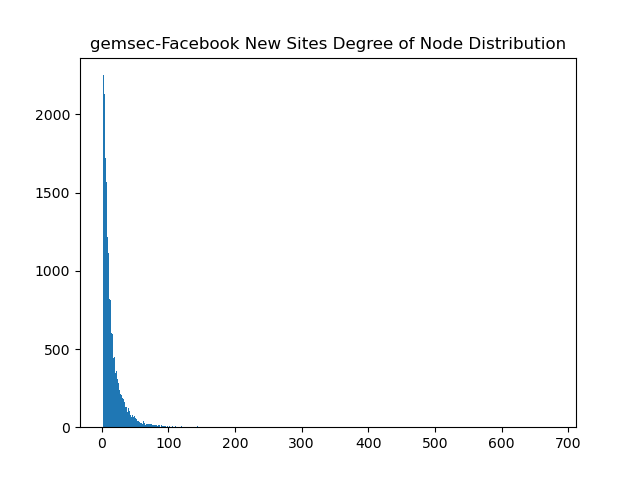}}\hfil 
\subfloat{\includegraphics[width=4cm]{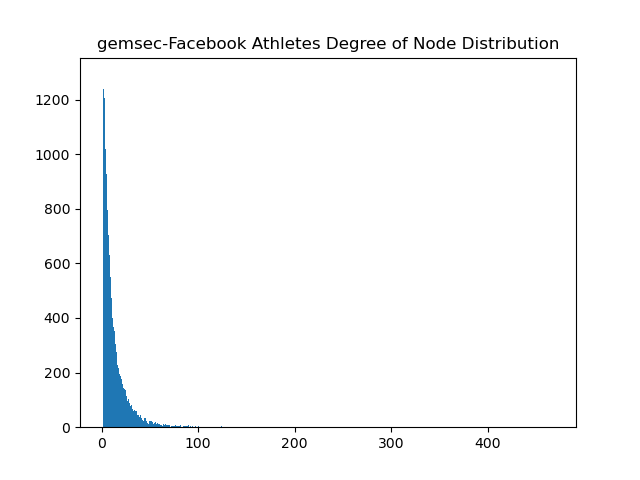}}\hfil
\subfloat{\includegraphics[width=4cm]{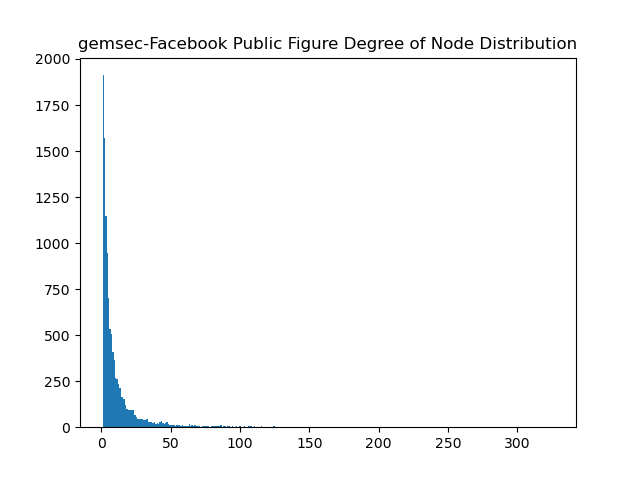}}\hfil 
\subfloat{\includegraphics[width=4cm]{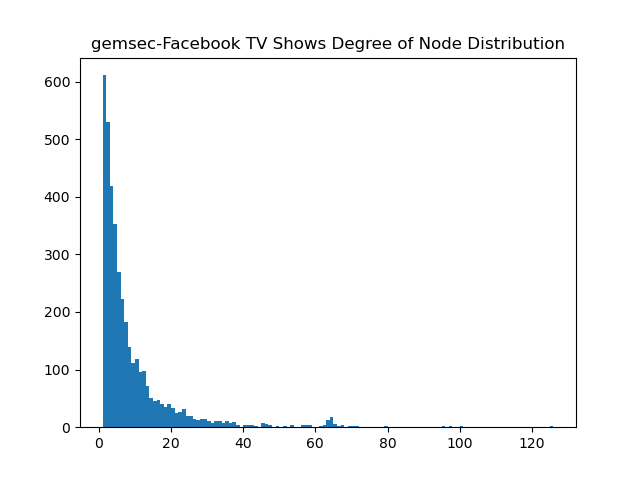}}\hfil
\subfloat{\includegraphics[width=4cm]{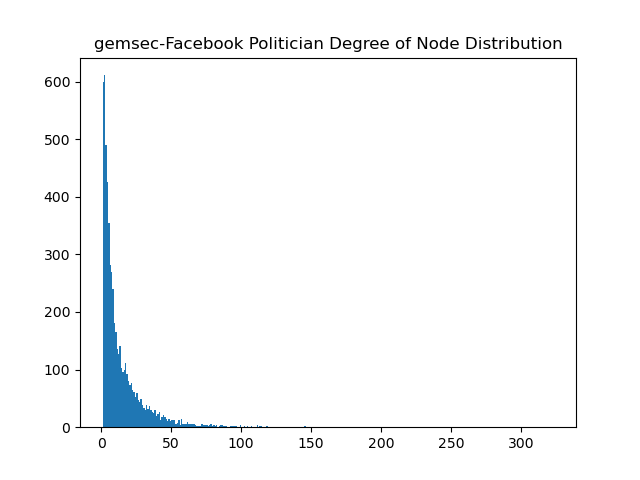}}\hfil
\subfloat{\includegraphics[width=4cm]{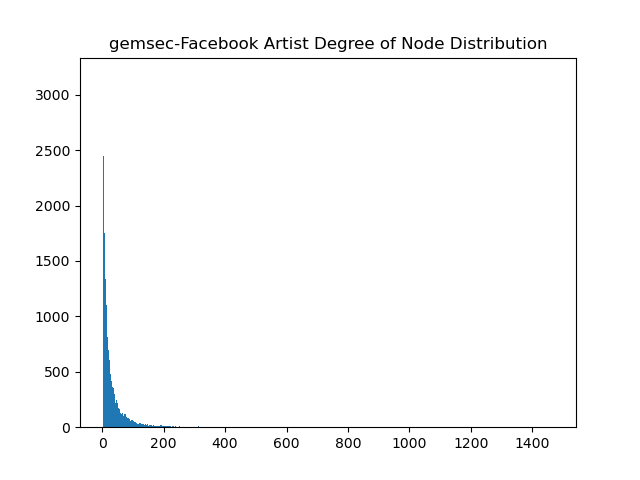}}\hfil
\subfloat{\includegraphics[width=4cm]{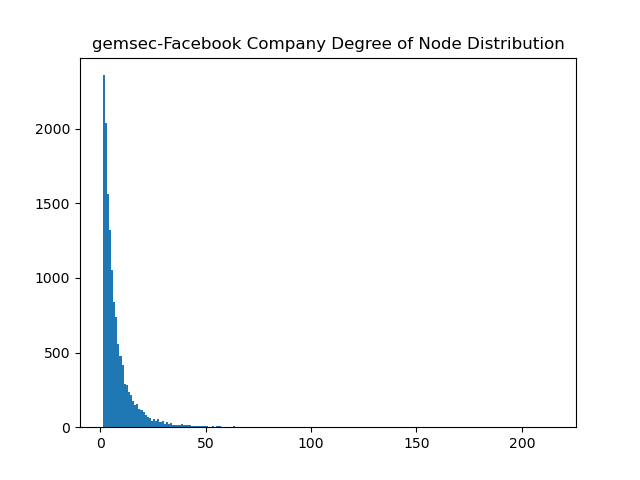}}
\caption{gamesec-Facebook}\label{fb}
\end{figure}
\begin{figure}
\centering
\subfloat{\includegraphics[width=3cm]{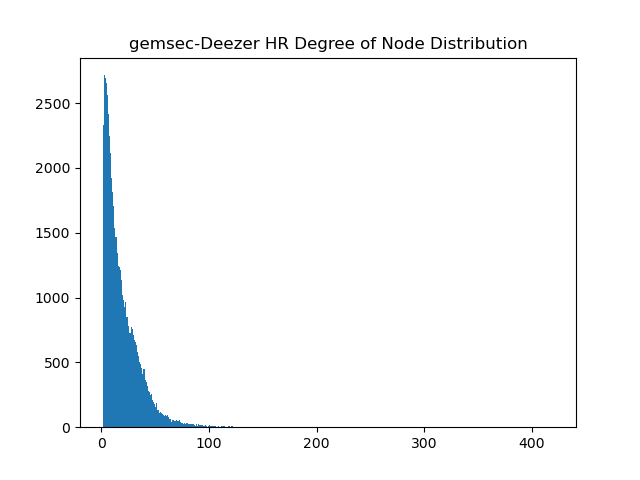}}
\subfloat{\includegraphics[width=3cm]{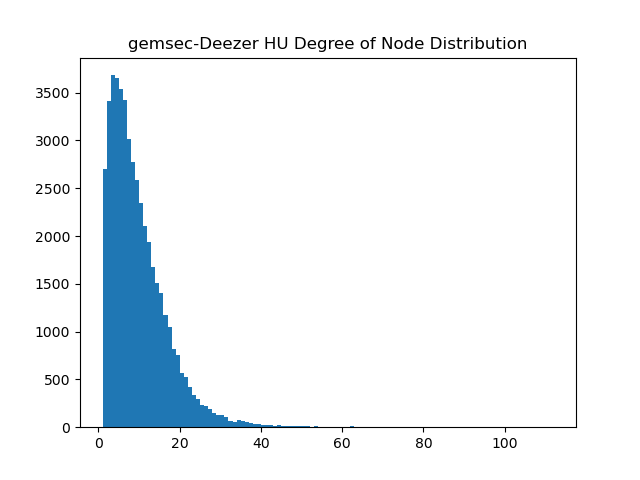}} 
\subfloat{\includegraphics[width=3cm]{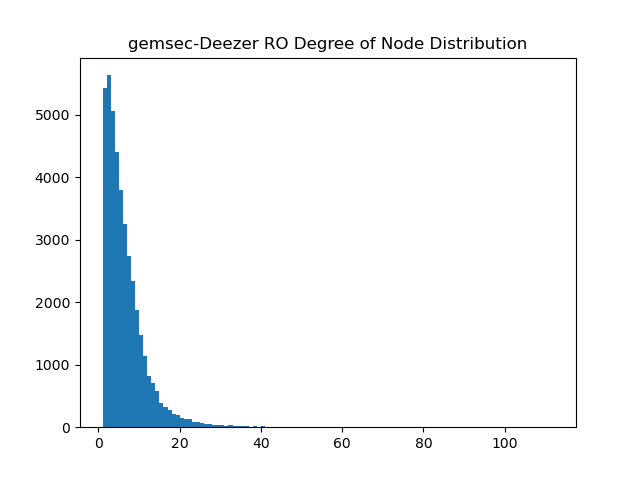}}
\caption{gamesec-Deezer}\label{deezer}
\end{figure}
\begin{figure}
\centering
\subfloat{\includegraphics[width=3cm]{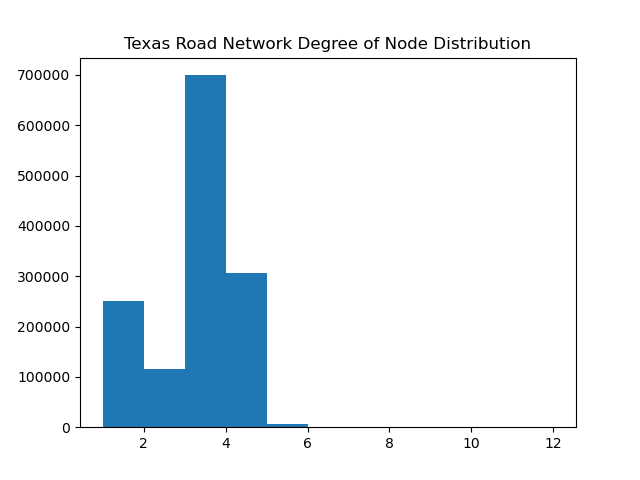}}
\subfloat{\includegraphics[width=3cm]{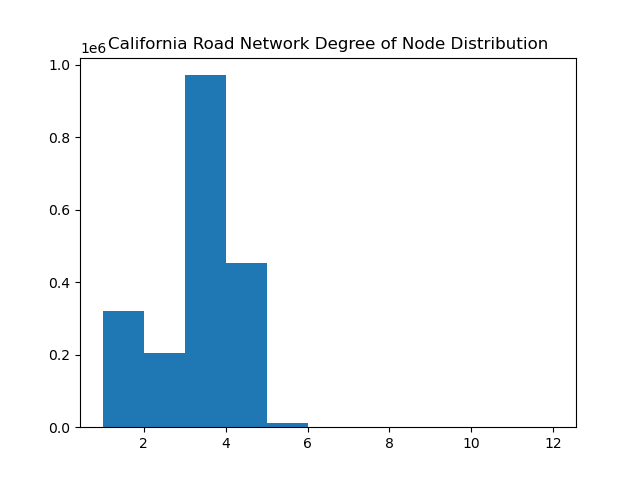}} 
\subfloat{\includegraphics[width=3cm]{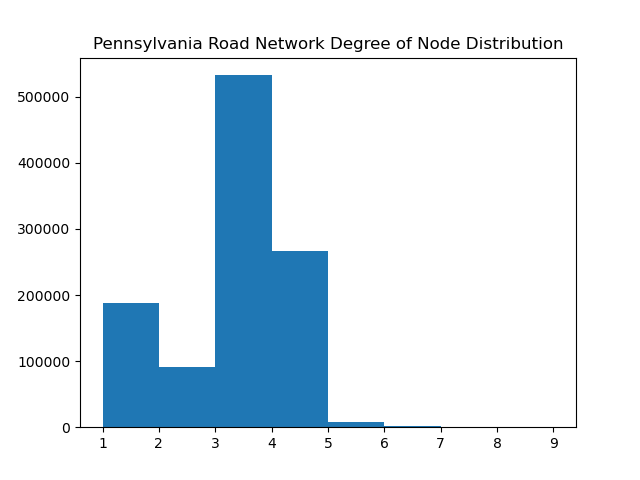}}
\caption{Road Network}\label{roadtable}
\end{figure}
\begin{figure}[H]
\centering
\includegraphics[width=12cm]{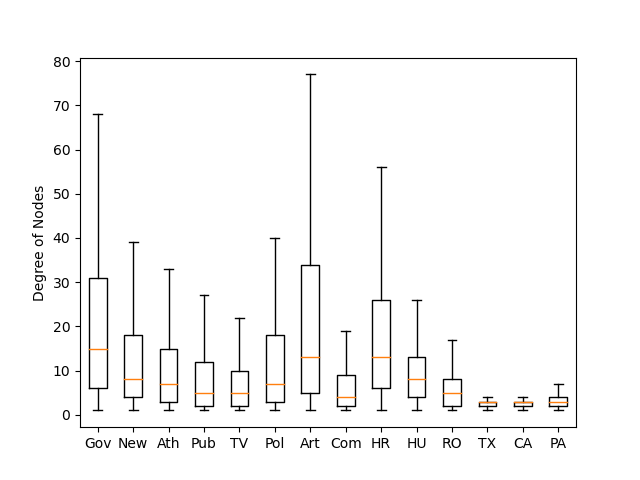}
\caption{Real-World Graphs Boxplot}
\label{fig:box}
\end{figure}
\noindent After testing with the deterministic structure graph and random graph generated by the computer, we can claim that algorithm 2 costs less time than algorithm 1 does, and the advantage can expand as the graph gets larger and average degree gets higher. Moreover, in a group of similar graphs, algorithm 2 works better when skewness and kurtosis of degree distribution are low since it prefers a high degree of node. 
\\\\
Next, we would test algorithm 2 on real-world problems to see if the above advantages continue. We select two social network data sets \cite{rozemberczki2019gemsec} and one road network data set \cite{road} from Stanford Large Data Network Collection.
\\\\
In figures \ref{fb}-\ref{roadtable}, are the degree distributions of three real-world data sets. we notice that the first two social networks data sets have high skewness and an extremely large degree of nodes that are hard to observe from graphs. However, road networks have a low degree in general. In tables \ref{table:fb}-\ref{table:road}, time ratio for gamesec-Facebook and gamesec-Deezer are relatively stable. While the nodes and edges vary widely, the average degree is high enough to ensure that algorithm 2 outperforms algorithm 1. When the average degree of graphs is low, such as in the road networks, the advantages of algorithm 2 are not as pronounced. This can be explained by the type of these data sets. The first two tables are social networks that have clusters and a high degree of nodes that can provide redundant edges for algorithm 2 to remove and save time. However, in road network graphs, algorithm 2 can only remove edges at cross-intersections and T-intersections. The property of the road network constrains the effectiveness of algorithm 2. Therefore, in real-world examples, algorithm 2 still works better on graphs that have a high average degree of nodes. In figure \ref{fig:box}, we pre-processed the data by removing the outliers to better observe the boxplots. The higher height of mean would potentially bring the time ratio in tables lower, and the ratio of IQR length and boxplot length seem to be stable throughout the graph which matches with the distribution graphs. Overall, we can conclude that algorithm 2 performs increasingly well on larger graphs with a higher average degree of nodes.

\section{Conclusion and Future Works}
Finding an optimal path cover in any graph is an NP-complete problem, yet there have been many works done to use approximation algorithms to estimate the optimal path cover. In this paper, we first introduce the $\frac{1}{2}$-Approximation Covering Algorithm as a fundamental path cover algorithm. A greedy algorithm that works with an undirected graph, it guarantees a path cover that obtains at least $\frac{1}{2}$ weight of the optimal path cover with time complexity of $O(M \cdot LogM)$. Since the steps of this algorithm are straightforward and deterministic, we optimize it by adding a step of removing redundant edges that connect to the existing path. We go over the proof in \cite{pathcover} to show that the $\frac{1}{2}$ theoretical bound still works for our optimized algorithm and visualize the process of removing redundant edges with a simple Watts-Strogatz graph. We then discuss our test results on both the deterministic version of Watts-Strogatz graphs, ring structured graphs, and Erdos-Renyi graphs, random graphs. From both the numerical results of time cost by algorithms and the visualizations, we observe that algorithm 2 works well on both types of graphs, especially as the average degree increases. Then we show that algorithm 2 is also effective on most real-world problems, except the low-degree sparse graphs such as the road networks. Therefore, a high degree of nodes in the graph would lead to a satisfying performance of algorithm 2.
\\\\
Overall, the sequential optimization algorithm attains superior performance than the primary algorithm with low computational time and stable time ratio, and the advantage would expand as the graph and its average degree gets larger. Because modern problems tend to have large data sets, the new algorithm will certainly be competitive if we bring it to applications, especially problems in the form of high degree networks, such as social activities and e-commerce. It is also applicable to graph theory problems, such as cocomparability graphs, cographs, interval graphs, block graphs, and permutation graphs. Each structure may require some changes of the algorithm, but it is a general approach in many cases.
\\\\
Even though the program thus far is good enough to use and meet the time complexity bound, we still want to explore if there are other potential data structures, functions, and programming tips of Python that we can implement to optimize our current program. Besides, we want to rewrite the current program in other languages, such as MATLAB or JAVA, which have better computation packages or work better with hardware directly to perform better. Additionally, we plan to take a step ahead to parallel optimization as we can separate one graph into multiple graphs by removing the local clustering and solve all the graphs at once to generate the path cover which would bring the efficiency to a higher level.
{\footnotesize  
\medskip
\medskip
\vspace*{1mm} 
 
\noindent {\it Junyuan Lin}\\  
Loyola Marymount University\\
1 LMU Drive\\
Los Angeles, CA 90045\\
E-mail: {\tt Junyuan.lin@lmu.edu}\\ \\  

\noindent {\it Guangpeng Ren}\\  
Loyola Marymount University\\
1 LMU Drive\\
Los Angeles, CA 90045\\
E-mail: {\tt gren@lion.lmu.edu}\\ \\}

\end{document}